\documentclass [a4paper,12pt]{article}
\usepackage{russ}
\usepackage[left=1cm,right=1cm,top=1cm,bottom=2cm]{geometry} 
\usepackage{amsmath}
\usepackage{amssymb}
\usepackage{amsthm}
\usepackage{amstext}
\usepackage{hyperref}
\usepackage{tikz}
\usetikzlibrary{calc}
\newtheorem{definition}{Определение}
\usepackage{float}
\usepackage{wrapfig}
\usepackage{mathptmx}       
\usepackage{helvet}         
\usepackage{courier}        

\title{Сравнение устойчивости процедур анализа характеристик сетевых моделей фондовых рынков.}
\author{Семенов Д.П., Калягин В.А., Колданов П.А., Бацын М.В., Голованова С.В.,\\ Воронина М.А.}
\begin{document}

\maketitle
\begin{abstract}Под сетевой моделью фондового рынка понимается полный взвешенный граф, вершины которого соответствуют доходностям рыночных активов, а веса ребер задаются мерой их зависимости. Для фильтрации ключевой информации из сетевой модели выделяются подграфы, которые мы будем называть сетевыми структурами. К популярным сетевым структурам в настоящее время относятся граф рынка, клики и независимые множества графа рынка, максимальное остовное дерево. При выделении таких структур и других характеристик сетевой модели из наблюдений неизбежно появляются ошибки, которые связаны со случайным характером данных, с конечностью времени наблюдения и с процедурами идентификации сетевых структур. Анализ таких ошибок в настоящей работе основан на вероятностной модели доходностей акций, которая описывается классом эллиптических распределений. Основное внимание сосредоточено на анализе устойчивости и сравнении характеристик двух типов процедур идентификации: широко известные процедуры, основанные на выборочном аналоге коэффициентов корреляции Пирсона и недавно предложенные процедуры, основанные на выборочных аналогах вероятностей совпадения знаков.
\end{abstract}
{\bf Ключевые слова:} фондовый рынок, сетевая модель, сетевые структуры, сеть случайных величин, процедуры идентификации, гистограмма распределения весов ребер, граф рынка, распределение степеней вершин графа рынка, клики и независимые множества графа рынка, топология максимального остовного дерева.
\section{Введение}
В последнее время всё большее распространение получают методы анализа фондового рынка, основанные на построении соответствующей сетевой модели\cite{Mantegna1999}, \cite{Boginski2005}, \cite{Tumminello2010}. Под сетевой моделью фондового рынка понимается полный взвешенный граф, вершины которого соответствуют доходностям рыночных активов (акций), а веса ребер задаются мерой их зависимости. Для фильтрации наиболее ценной информации из такой модели выделяются подграфы, которые называются сетевыми структурами. К популярным сетевым структурам можно отнести граф рынка, клики и независимые множества графа рынка \cite{Boginski2005} и максимальное остовное дерево \cite{Mantegna1999}. К настоящему времени известно достаточно большое количество работ по применению такого подхода к анализу рынков различных стран и интерпретации полученных результатов \cite{Coronnello2005}, \cite{GarasArgyrakis2007}, \cite{Huang2009}, \cite{Jung2006}, \cite{Tabak2010}, \cite{Tee2005}, \cite{Wang2012}, \cite{Tumminello2005},\cite{Визгунов2012}, \cite{Koldanov2014}, \cite{Bonanno2001}, \cite{Boginski2003}, \cite{Boginski2006}. Вместе с тем, в этих исследованиях отсутствует анализ достоверности полученных выводов. Такой анализ необходим так как наблюдения над акциями представляют собой наблюдения над случайными величинами \cite{Ширяев1998}, \cite{Фантаццини2011}. Следовательно, при выделении сетевых структур и анализе других характеристик сетевой модели по наблюдениям неизбежно появляются ошибки, которые связаны со случайным характером данных, с конечностью времени наблюдения и с применяемыми процедурами идентификации сетевых структур. Первой работой, в которой задача идентификации сетевых структур рассматривалась с учетом случайного характера наблюдений является работа \cite{Koldanov2013}. Настоящая работа является продолжением исследований в этом направлении.
\\
Математической основой проводимых исследований является понятие сети случайных величин \cite{ЖНЭА2017},\cite{Kalyagin2017}. Сетью случайных величин называется пара $({X, \gamma})$, где $X=(X_1,X_2,...,X_N)$ — случайный вектор, $\gamma=\gamma(X,Y)$ — мера зависимости случайных величин $X,Y$. В сетевой модели, порожденной сетью случайных величин, вес ребра $(i,j)$ задается значением $\gamma(X_i,X_j)$. Для построения сети случайных величин фондового рынка и порожденной сетевой модели, вектор $X=(X_1,X_2,...,X_N)$ задаётся совместным распределением доходностей рыночных активов. Для описания такого распределения в настоящее время используются многомерные эллиптические распределения \cite{Gupta2013} и распределения, заданные различными копула-функциями \cite{Фантаццини2011}. В широкий класс эллиптических распределений входят, в частности, многомерные нормальные распределения и распределения Стьюдента с тяжелыми хвостами. В качестве меры зависимости доходностей обычно используется коэффициент корреляции Пирсона. Вместе с тем, в работах \cite{Bautin2013a}, \cite{Bautin2014} было показано, что неопределенность процедур идентификации графа рынка и максимального остовного дерева, основанных на выборочных корреляциях Пирсона, может существенно возрастать при отклонении совместного распределения доходностей от нормального. В работе \cite{Bautin2013b} предложена новая мера зависимости доходностей, основанная на вероятности совпадения знаков. В работах \cite{Bautin2013a}, \cite{Bautin2014} экспериментально установлено, что процедуры идентификации графа рынка и максимального остовного дерева, основанные на вероятности совпадения знаков, обладают свойством устойчивости при отклонении распределения от нормального. В работе \cite{Kalyagin2017} теоретически доказана устойчивость статистических процедур идентификации графа рынка в сетевых моделях, порожденных сетями случайных величин в классе эллиптических распределений. В \cite{ЖНЭА2017},\cite{Kalyagin2017} доказано, что этот вывод справедлив и для процедур идентификации максимального остовного дерева. Вместе с тем, остаётся открытым вопрос насколько сильно отличаются характеристики устойчивости процедур идентификации клик и независимых множеств, а также процедур оценивания других характеристик сетевых моделей и сетевых структур.\\
Целью настоящей работы является сравнение характеристик устойчивости двух типов процедур идентификации: широко известные процедуры, основанные на выборочном коэффициенте корреляций Пирсона и процедуры, основанные на оценках вероятностей совпадения знаков в классе эллиптических распределений. При этом изучаются такие сетевые структуры как: максимальные клики и максимальные независимые множества графа рынка и такие характеристики как: распределение весов рёбер сетевой модели (полного взвешенного графа) фондового рынка, распределение степеней вершин в графе рынка, топология степеней вершин максимального остовного дерева.\\
Работа организована следующим образом: в разделе \ref{sec2} описана вероятностная модель, сеть случайных величин, введены основные определения и обозначения; в разделе \ref{sec3} приведены формулировки задач; в разделе \ref{sec4} изложен общий подход к сравнению истинных и выборочных характеристик сетевых моделей и их структур; в разделе \ref{sec5} приведены результаты статистического моделирования, основанные на анализе рынков 8 стран и обсуждаются полученные результаты.

\section{Основные определения и обозначения} \label{sec2}
\begin{definition}\label{random_variables_network_model}
Сетью случайных величин будем называть пару $(X,\gamma)$, где $X=(X_1,\ldots,X_N)$-вектор случайных величин, а $\gamma=\{\gamma_{i,j}:i,j=1,\ldots,N;i\neq j\}$ мера зависимости между случайными величинами $X_i,X_j$.
\end{definition}

В зависимости от распределения вектора $X$ и выбранной меры зависимости $\gamma$ можно рассматривать различные сети случайных величин.
Будем предполагать, что вектор $X$ имеет эллиптическое распределение, плотность которого имеет вид \cite{Anderson2003}:
\begin{equation} \label{ell_con}
f(x)=|\Lambda|^{-1}g((x-\mu)'\Lambda^{-1}(x-\mu)),
\end{equation}
 где $\Lambda$ - положительно определенная матрица, $g(x)\geq0$, а функция распределения $g(x'x)$ удовлетворяет условию нормировки $$\int_{-\infty}^{\infty}g(x'x)d(x)=1$$
К такому классу распределений относятся: многомерное нормальное $$f(x)=\frac{1}{(2\pi)^{n/2}|\Lambda|^{1/2}}exp^{\frac{-1}{2}(x-\mu)'\Lambda^{-1}(x-\mu)}$$ многомерное распределение Стьюдента $$t_{\nu}(x)=\frac{1}{\sqrt{(\nu\pi)^N*det\Lambda}}\frac{Г(\frac{\nu+N}{2})}{Г(\frac{\nu}{2})}(1+\frac{(x-\mu)'\Lambda^{-1}(x-\mu)}{\nu})^{\frac{\nu+N}{2}}$$ и их смесь при одном и том же векторе $\mu$ и матрице $\Lambda$, где $\nu$ - количество степеней свободы. \\
В настоящей работе рассматриваются следующие сети случайных величин.
\begin{definition}\label{random_variables_sign_network_pearson_measure}
Эллиптической сетью корреляций Пирсона будем называть сеть случайных величин $(X,\gamma^{P})$, в которой вектор $X=(X_1,\ldots,X_N)$ имеет распределение с плотностью, постоянной на многомерных эллипсоидах, а мера зависимости $\gamma_{i,j}^{P}$ между случайными величинами $X_i,X_j$ задается коэффициентом корреляции Пирсона $$\rho_{i,j}=\frac{E(X_i-EX_i)(X_j-EX_j)}{\sqrt{DX_iDX_j}}$$
\end{definition}

Такие сети (без предположения о распределении $X$) широко применяются в задачах генетики и финансового рынка \cite{Mantegna1999}, \cite{Boginski2003}, \cite{Boginski2005}, \cite{Boginski2006}.

\begin{definition}\label{random_variables_sign_network}
Эллиптической сетью вероятностей совпадения знаков или знаковой сетью будем называть сеть случайных величин $(X,\gamma^{Sg})$, в которой вектор $X=(X_1,\ldots,X_N)$ имеет распределение с плотностью, постоянной на многомерных эллипсоидах, а мера зависимости $\gamma_{i,j}^{Sg}$ между случайными величинами $X_i,X_j$ задается вероятностью совпадения знаков $$p^{i,j}=P((X_i-\mu_i)(X_j-\mu_j)>0)$$
\end{definition}

Сети случайных величин порождают сетевые модели, которые представляют собой простой полный неориентированный взвешенный граф  $G=(V, E, \gamma)$, где $V=\{1,2,\ldots,N\}$ - множество вершин, которые описываются случайными величинами $X_1,X_2,\ldots,X_N$, E - множество ребер с весами, заданными мерой $\gamma$.

Изучение сетевых моделей $G=(V,E, \gamma)$ естественно свести к изучению ключевых характеристик соответствующих графов. В теории графов предложено достаточно большое количество таких характеристик: отсеченный граф, клики, независимые множества, максимальное остовное дерево, степени вершин, центральность, диаметр и т.д. \\
В настоящей работе исследуются характеристики графов, удовлетворяющие следующим определениям:
\begin{definition} \label{hist}
Распределением весов ребер называется функция $h(x)=m$, где $m$ - число ребер веса $x$.
\end{definition}
\begin{definition} \label{def_threshold_graph}
Отсеченным графом (графом рынка, MG) сетевой модели $G=(V,E,\gamma)$ называется подграф $G'(\gamma_0)=(V',E'):V'=V; E'\subseteq E, E'=\{(i,j):\gamma_{i,j}> \gamma_0\}$, где $\gamma_0$ - некоторый порог.
\end{definition}
Подчеркнем, что граф рынка представляет собой простой неориентированный граф без весов и без петель. Вместе с графом рынка часто изучаются его клики и независимые множества.
\begin{definition}\label{Vdeg}
Под распределением степеней вершин графа рынка понимается таблица $2\times N$, где в первой строке указаны возможные значения степеней вершин $0, 1,\ldots,N-1$, а во второй строке указано число вершин $\nu_i$ степени $i, i=0,\ldots,N-1$.
\end{definition}
\begin{definition}\label{clique}
Кликой графа рынка $G=(V,E)$ называется полный подграф графа $G$, т.е. подграф $G'=(V',E'):V'\subset V, E'\subset E:\forall i, j\in V'\Rightarrow (i,j)\in E'$
\end{definition}
\begin{definition}\label{maximum_clique}
Клика $G_1=(V_1 ,E_1)$ называется максимальной(MC) (по размеру), если для любой другой клики $G_2=(V_2 ,E_2)$ выполняется: $|V_1|\geq|V_2|$
\end{definition}

\begin{definition}\label{independent_set}
Независимым множеством (IS) графа рынка $G=(V,E)$ называется пустой подграф графа $G$, т.е. подграф $G_1=(V_1 ,E_1):V_1\subset V, E_1\subset E:\forall i, j\in V_1\Rightarrow (i,j)\notin E_1$
\end{definition}
\begin{definition}\label{maximum_independent_set}
Независимое множество $G_1=(V_1 ,E_1)$ называется максимальным (MIS) (по размеру), если для любого другого независимого множества $G_2=(V_2 ,E_2)$ графа $G$ выполняется: $|V_1|\geq|V_2|$.
\end{definition}
Семейство $\{MG(\gamma_0):\gamma_0\in R^1\}$ содержит наиболее полную информацию о сетевой модели, в частности, о сетевой модели рынка. При этом клики и независимые множества характеризуют кластерную структуру рынка. Кроме того, размер максимальных клик является показателем глобализации, а размер максимального независимого множества является показателем степени 'свободы' на рынке.
\begin{definition} \label{def_MST}
Максимальным остовным деревом (MST) сетевой модели $G=(V,E,\gamma)$ называется дерево (граф без циклов) $G'=(V',E'):V'=V; E'\subset E; |E'|=|V|-1;$ такое, что $\sum_{(i,j)\in E'}\gamma_{i,j}$ максимальна.
\end{definition}
\begin{definition} \label{def_MST_topology}
Топологией MST будем называть последовательность степеней вершин MST, упорядоченную в возрастающем порядке.
\end{definition}

\section{Задачи идентификации характеристик сетевых моделей фондовых рынков}\label{sec3}
При практическом построении характеристик сетевых моделей распределение вектора $X$ и значение $\gamma_{i,j}$ неизвестны. Доступными данными являются наблюдения за доходностями рыночных активов. Под проблемой идентификации характеристик сетевых моделей в настоящей работе понимается задачи построения таких характеристик по наблюдениям. В качестве модели наблюдений используется повторная выборка $x_i(t), i=1,\ldots,N; t=1,\ldots,n$ конечного объема из распределения случайного вектора $X=(X_1,\ldots,X_N)$. \\
\begin{enumerate}

\item Задача оценки распределения весов ребер $h(x)$ заключается в построении гистограммы оценок весов рёбер.
\item Задача идентификации графа рынка заключается в выборе одной из гипотез:
\begin{equation}\label{problem_formulation_threshold_graph_example}
\begin{array}{ll}
&H_{1}^{TG}:\gamma_{i, j}\leq\gamma_{0},\forall (i, j), \ i < j,\\
&H_{2}^{TG}:\gamma_{1,2}>\gamma_{0},\gamma_{i, j}\leq\gamma_{0},\forall(i, j)\neq(1,2), \ i < j, \\
&H_{3}^{TG}:\gamma_{1,2}>\gamma_{0},\gamma_{13}>\gamma_{0},\gamma_{i, j}\leq \gamma_{0},\forall(i, j)\neq(1,2),(i, j)\neq(1,3),\\
&\ldots \\
&H_{L}^{TG}:\gamma_{i, j}>\gamma_{0},\forall(i, j), \ i < j,
\end{array}
\end{equation}

При этом гипотеза $H_1^{TG}$ соответствует пустому графу рынка $G'(\gamma_0)$, гипотеза $H_2^{TG}$ - графу рынка $G'(\gamma_0)$ с одним ребром $(1,2)$, и т.д., гипотеза $H_L^{TG}$ соответствует полному графу рынка $G'(\gamma_0)$. Для сетевой модели на $N$ вершинах число гипотез равно $L=2^{\frac{N(N-1)}{2}}$
\item Задача оценки степеней вершин в графе рынка.
\item Задача идентификации максимальных клик в графе рынка.
\item Задача идентификации максимальных независимых множеств в графе рынка.
\item Задача идентификации MST.
\\
Пусть $E_1=\{(i,j):(i,j)\in E\}$ - подмножество ребер, образующих  остовное дерево, $E_2=\{(i,j):(i,j)\in E\}$ - другое подмножество ребер, образующих  остовное дерево и т.д.  Обозначим множество всех таких подмножеств $E_{MST}=\{E_1,E_2,\ldots,E_{L_{MST}}\}$.
В соответствии с определением \ref{def_MST} задача идентификации максимального остовного дерева (при условии, что оно единственно) по наблюдениям $x_i(t);i=1,\ldots,N, t=1,\ldots,n$ может быть сформулирована как задача выбора одной из многих статистических гипотез

\begin{equation}\label{problem_formulation_MST_example}
\begin{array}{ll}
&H_{1}^{MST}:\sum_{(i,j)\in E_1}\gamma_{i, j}>\sum_{(i,j)\in E_k}\gamma_{i, j}:\forall E_k\in E_{MST},E_k\neq E_1\\
&H_{2}^{MST}:\sum_{(i,j)\in E_2}\gamma_{i, j}>\sum_{(i,j)\in E_k}\gamma_{i, j}:\forall E_k\in E_{MST},E_k\neq E_2\\
&H_{3}^{MST}:\sum_{(i,j)\in E_3}\gamma_{i, j}>\sum_{(i,j)\in E_k}\gamma_{i, j}:\forall E_k\in E_{MST},E_k\neq E_3\\
&\ldots \\
&H_{L_{MST}}^{MST}:\sum_{(i,j)\in E_{L_{MST}}}\gamma_{i, j}>\sum_{(i,j)\in E_k}\gamma_{i, j}:\forall E_k\in E_{MST},E_k\neq E_{L_{MST}}
\end{array}
\end{equation}
Согласно формуле Кэли число остовных деревьев в полном графе на $N$ вершинах, и следовательно, число различаемых гипотез равно $L_{MST}=N^{N-2}$.
\item Задача оценки топологии MST.
\end{enumerate}
\section{Меры различия истинных и выборочных характеристик}\label{sec4}
Пусть $\gamma_{i,j}, i,j=1,\ldots,N$ - истинное значение меры зависимости между случайными величинами $X_i$ и $X_j$ (вес ребра между вершинами $i$ и $j$ в сетевой модели). Сетевую модель, построенную на основе $\gamma_{i,j}, i,j=1,\ldots,N$ будем называть {\it истинной} сетевой моделью. Характеристики этой сетевой модели, определенные в разделе 2 будем называть {\it истинными} характеристиками сетевой модели.\\
Пусть $\hat{\gamma}_{i,j}, i,j=1,\ldots,N$ - оценка значения меры зависимости между случайными величинами $X_i$ и $X_j$, построенная по выборке $x_i(t), i=1,\ldots,N; t=1,\ldots,n$. Сетевую модель, построенную на основе $\hat{\gamma}_{i,j}, i,j=1,\ldots,N$ будем называть {\it выборочной} сетевой моделью. Характеристики этой сетевой модели будем называть {\it выборочными} характеристиками сетевой модели.\\
Введем меры различия между истинными и выборочными характеристиками соответствующих сетевых моделей.
\begin{definition} \label{def_difference_measure_of_histograms_of_weights_of_edges}
Под мерой различия истинного распределения весов ребер $h(x)$ и оценкой этого распределения (гистограммой $\hat{h}(x)$) будем понимать $E(|S-\hat{S}|)$, где $E$ - математическое ожидание, $S$ - площадь под кривой $h(x)$, $\hat{S}$ - площадь под кривой $\hat{h}(x)$.
\end{definition}
\begin{definition} \label{def_difference_measure_of_degrees_distributions}
Под мерой различия оценки распределения степеней вершин графа рынка и его истинного значения будем понимать $E(\sum_{i=0}^{N-1} |k_i -\hat{k}_i|)$, где $N$ – число вершин, $k_i$ – истинное число вершин степени $i$, $\hat{k}_i$ – ее оценка.
\end{definition}
\begin{definition} \label{def_difference_measure_of_cliques}
Под мерой различия максимальных выборочных и истинных клик (независимых множеств) будем понимать $E(\sum_{i=1}^N |v_i -\hat{v}_i|)$ математическое ожидание мощности симметрической разности множества вершин истинной клики и выборочной клики, где $N$ – число вершин, $v_i$ – индикатор, равный $1$, если вершина $i$ есть в клике (независимом множестве), а $\hat{v}_i$ – индикатор в выборочной клике (независимом множестве).
\end{definition}
\begin{definition} \label{def_difference_measure_of_MSTs}
Под мерой различия топологии выборочного и истинного MST будем понимать вероятность правильного определения топологии степеней вершин.
\end{definition}
\begin{definition} \label{stability}
Под устойчивостью введенных мер различия к изменению вероятностной модели распределения вектора $X$ будем понимать независимость этих мер от функции $g$.
\end{definition}

\section{Результаты исследования устойчивости процедур оценки характеристик сетевых моделей}\label{sec5}
В настоящем разделе приведены результаты исследования устойчивости двух типов процедур идентификации характеристик сетевых моделей. Исследования проводились в рамках эллиптической сети корреляции Пирсона. В качестве истинных сетевых моделей были использованы сетевые модели, построенные по наблюдениям за доходностями акций реальных рынков. Были проанализированы доходности $50$ наиболее доходных акций рынков России, Бразилии, Индии, Китая, Франции, Германии, США и Великобритании за 2010-2012 года. В качестве истинных значений весов ребер использовались значения коэффициентов корреляции Пирсона, построенные по наблюдениям за $1$ календарный год. В таблице \ref{table1} для примера приведен фрагмент матрицы $10\times10$ весов истинной сетевой модели рынка Великобритании за 2010 год.
\begin{table}[H]
\centering
\caption{Фрагмент истинной матрицы корреляций доходностей. Великобритания, 2010 год. Сеть Пирсона.}
\begin{tabular}{llllllllll}
1,00  & 0,12  & 0,00  & 0,10  & -0,05 & -0,14 & 0,04  & -0,02 & -0,02 & 0,22  \\
0,12  & 1,00  & 0,08  & -0,03 & -0,01 & -0,03 & -0,05 & -0,03 & 0,05  & -0,10 \\
0,00  & 0,08  & 1,00  & 0,04  & -0,06 & 0,02  & 0,02  & -0,04 & -0,04 & -0,03 \\
0,10  & -0,03 & 0,04  & 1,00  & -0,07 & 0,06  & 0,10  & 0,06  & 0,04  & 0,09  \\
-0,05 & -0,01 & -0,06 & -0,07 & 1,00  & 0,49  & 0,14  & 0,44  & 0,35  & 0,01  \\
-0,14 & -0,03 & 0,02  & 0,06  & 0,49  & 1,00  & 0,24  & 0,48  & 0,42  & -0,09 \\
0,04  & -0,05 & 0,02  & 0,10  & 0,14  & 0,24  & 1,00  & 0,30  & 0,15  & 0,00  \\
-0,02 & -0,03 & -0,04 & 0,06  & 0,44  & 0,48  & 0,30  & 1,00  & 0,45  & 0,04  \\
-0,02 & 0,05  & -0,04 & 0,04  & 0,35  & 0,42  & 0,15  & 0,45  & 1,00  & 0,01  \\
0,22  & -0,10 & -0,03 & 0,09  & 0,01  & -0,09 & 0,00  & 0,04  & 0,01  & 1,00
\end{tabular}
\label{table1}
\end{table}
На основе таких сетевых моделей за каждый год наблюдений для каждой страны были построены истинные характеристики. В общей сложности для 8 стран, 3 лет наблюдений и 5 характеристик было построено $8*3*5=120$ истинных характеристик сетевых моделей. Примеры истинных характеристик сетевых моделей приведены ниже.
\begin{figure}[H]
\label{fig1}
\center{\includegraphics*[scale=0.3]{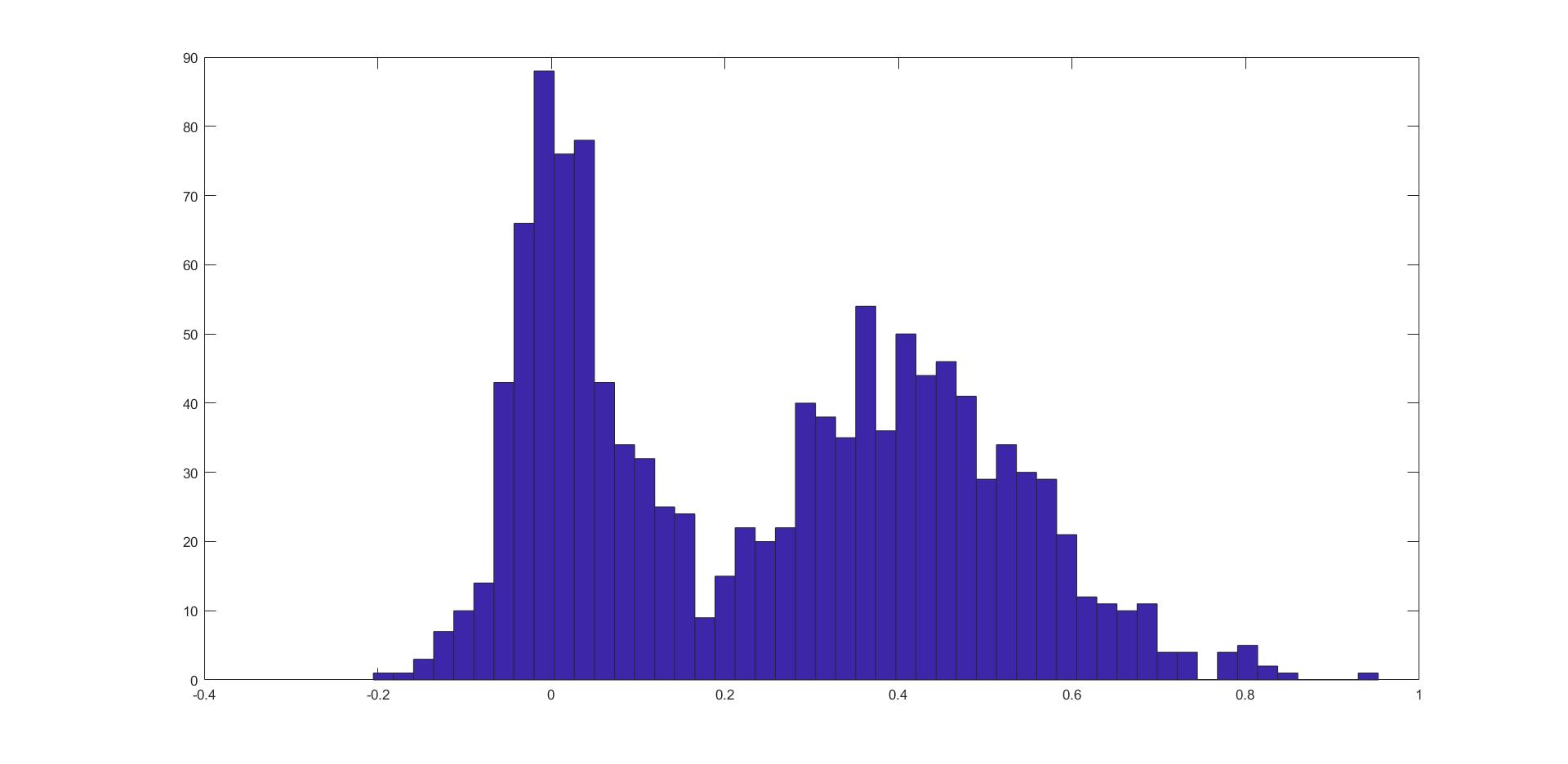}}
\caption{Истинная гистограмма весов ребер. Великобритания, 2010 год. Сеть Пирсона.}
\end{figure}
\begin{figure}[H]
\label{fig2}
\center{\includegraphics*[scale=0.3]{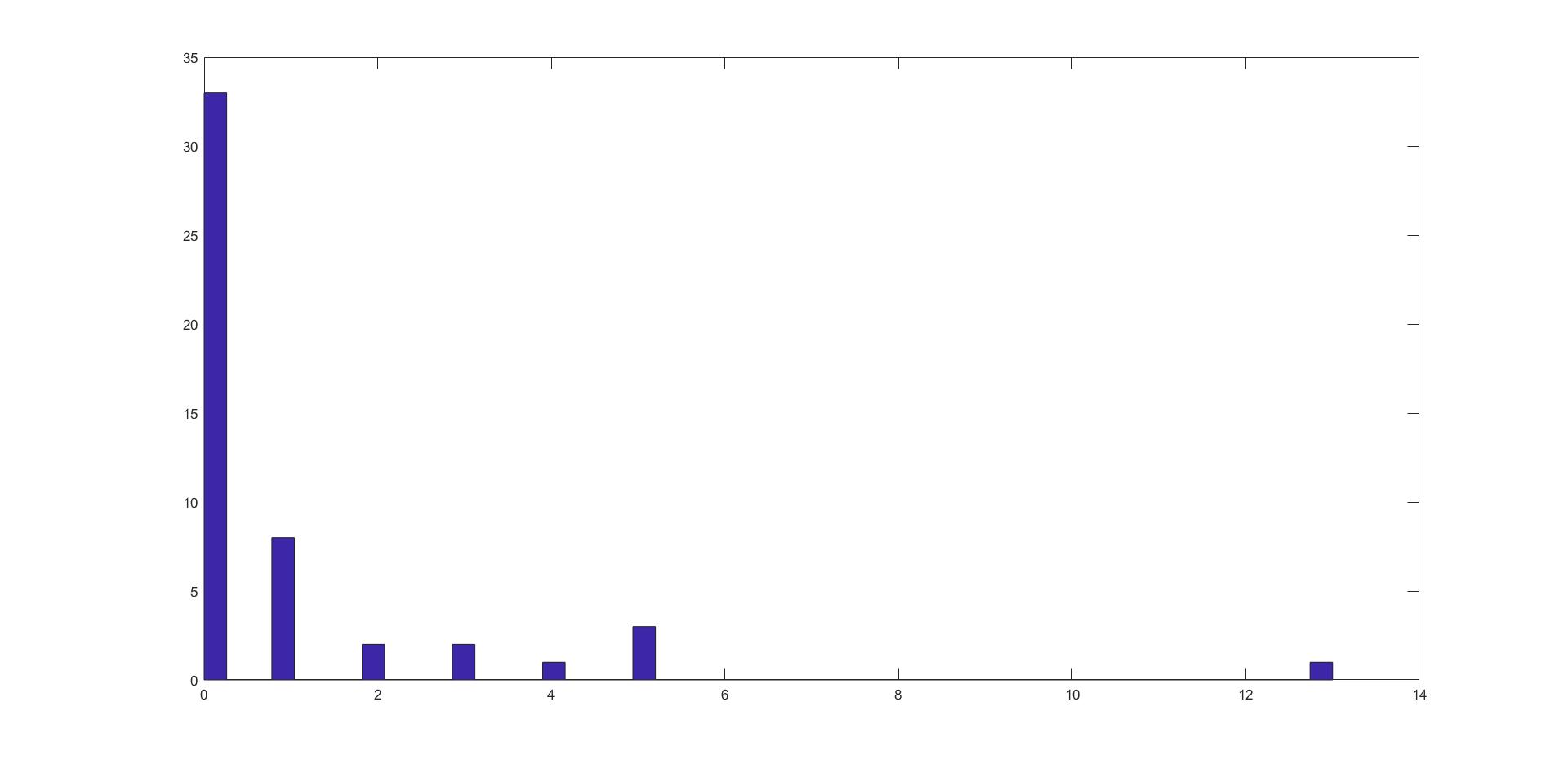}}
\caption{Истинное распределение степеней вершин в графе рынка при пороге 0.3. Великобритания, 2010 год. Сеть Пирсона.}
\end{figure}
\begin{figure}[H]
\label{fig3}
\center{\includegraphics*[scale=0.5]{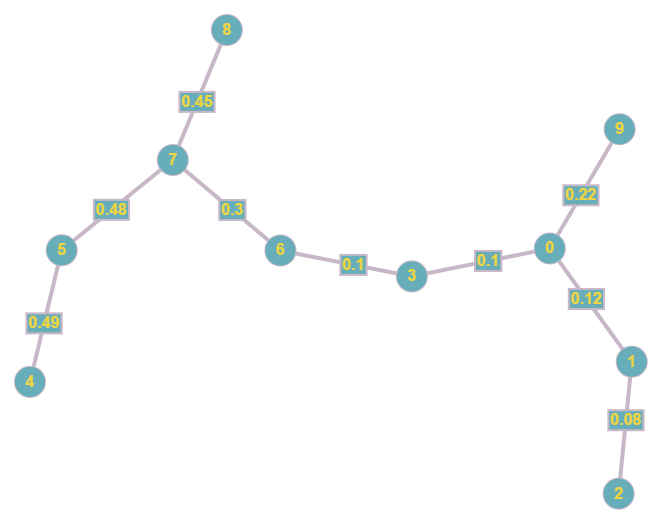}}
\caption{Истинное MST. Великобритания, 2010 год. Сеть Пирсона.}
\end{figure}
Анализ устойчивости рассматриваемых процедур идентификации характеристик сетевых моделей проводился методом статистического моделирования. Были сгенерированы многомерные выборки различного объема из распределения смеси с плотностью: $$f(x_1,...,x_N)=\gamma*f_{gauss}(x_1,...,x_N) + (1-\gamma)*f_{St,k}(x_1,...,x_N)$$ где $f_{gauss}(x_1,...,x_N)$ - N-мерное нормальное распределение, $f_{St, k}(x_1,...,x_N)$ - N-мерное распределение Стьюдента с $\nu=3$ степенями свободы. В качестве $\Lambda$ использовалась истинная матрица весов ребер. Эксперимент повторялся 10000 раз и усреднением находилась оценка введенных выше мер различия. Результаты проведения экспериментов представлены в виде кривых зависимости меры различия от $\gamma (\gamma=0, 0.1, 0.2, ... , 1)$  для каждого рынка и каждой сети.


\subsection{Распределение весов ребер}
Для оценки выборочных гистограмм использовались традиционные процедуры: выборочная корреляция Пирсона и частота совпадения знаков. Эксперимент заключался в следующем: $n=100$  раз генерировался N-мерный $(N=50)$ случайный вектор с распределением смеси с заданным $\gamma$. По наблюдениям вычислялась гистограмма и мера различия. Эксперимент повторялся 10000 раз и усреднением находилась оценка меры различия. Результаты проведения экспериментов представлены в виде кривых зависимости меры различия от $\gamma (\gamma=0, 0.1, 0.2, ... , 1)$  для каждого рынка и каждой сети.
\\
Результаты экспериментов показывают, что оценка распределения весов ребер, основанная на оценке вероятности совпадения знаков устойчива к изменению параметра смеси $\gamma$. При использовании для оценок выборочных коэффициентов корреляции Пирсона мера различия неустойчива к изменению параметра смеси $\gamma$. Аналогичные результаты, показывающие устойчивость процедур, основанных на оценках вероятности совпадения знаков справедливы и для рынков других стран и годов наблюдения.
\begin{figure}[H]
\label{fig4}
\center{\includegraphics*[scale=0.8]{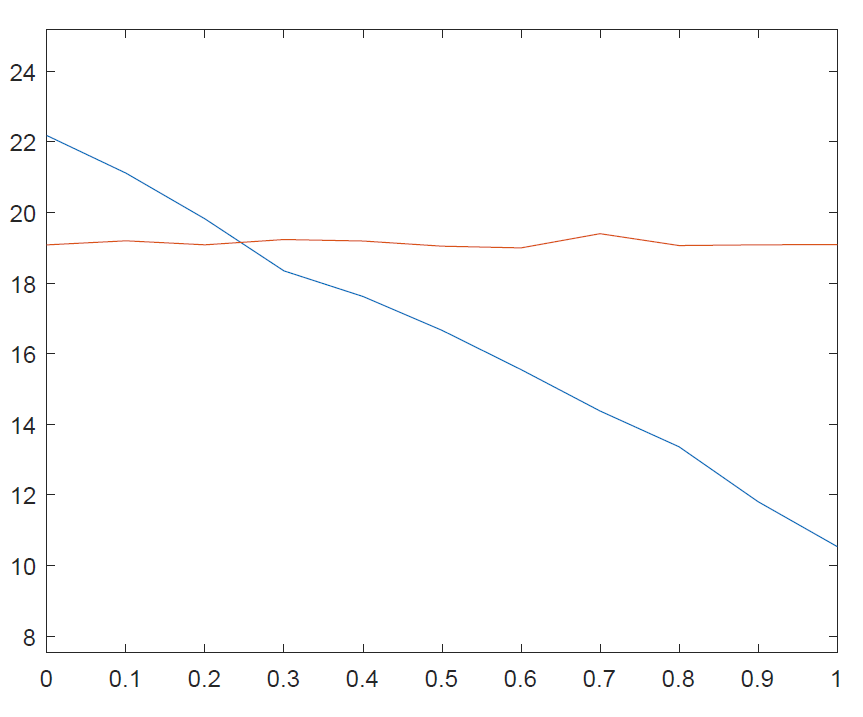}}
\caption{Зависимость меры различия \ref{def_difference_measure_of_histograms_of_weights_of_edges} от $\gamma$. Гистограмма весов ребер. Россия, 2012 год. Красным - сеть вероятности совпадения знаков, синим - сеть Пирсона.}
\end{figure}
\begin{figure}[H]
\label{fig5}
\center{\includegraphics*[scale=0.8]{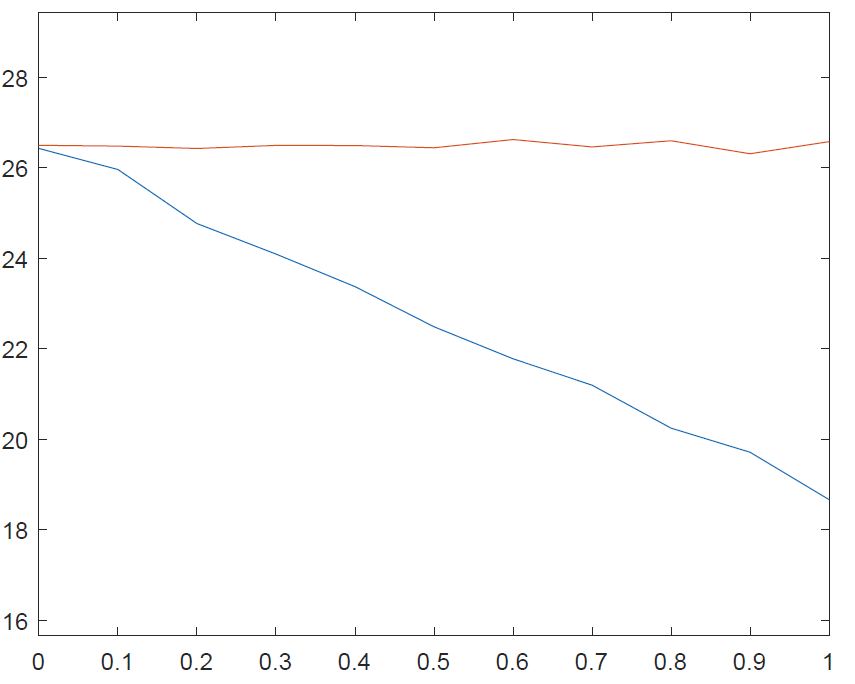}}
\caption{Зависимость меры различия \ref{def_difference_measure_of_histograms_of_weights_of_edges} от $\gamma$. Гистограмма весов ребер. Китай, 2010 год. Красным - сеть вероятности совпадения знаков, синим - сеть Пирсона.}
\end{figure}
\subsection{Распределение степеней вершин}
Эксперимент заключался в следующем: $n (n=100, 250)$  раз генерировался N-мерный $(N=50)$ случайный вектор с распределением смеси с заданным $\gamma$. По наблюдениям вычислялись выборочные корреляции Пирсона и частоты совпадения знаков, строился граф рынка \ref{def_threshold_graph} с порогом $\gamma_0 (\gamma_0=0,1; 0,3; 0,5)$, находилось распределение степеней вершин и мера различия. Эксперимент повторялся 10000 раз и усреднением находилась оценка меры различия. Результаты проведения экспериментов представлены в виде кривых зависимости меры различия от $\gamma (\gamma=0, 0.1, 0.2, ... , 1)$ для каждого рынка, каждой сети и каждого значения количества наблюдений.
\\
Как и в случае с распределением весов ребер, результаты экспериментов показывают, что вероятность совпадения знаков оказывается устойчивой к изменению параметра смеси $\gamma$, а оценки выборочных коэффициентов корреляции Пирсона - нет. Также стоит отметить, что значение меры различия, как правило, быстро уменьшается с ростом числа наблюдений.
\begin{figure}[H]
\label{fig6}
\center{\includegraphics*[scale=0.7]{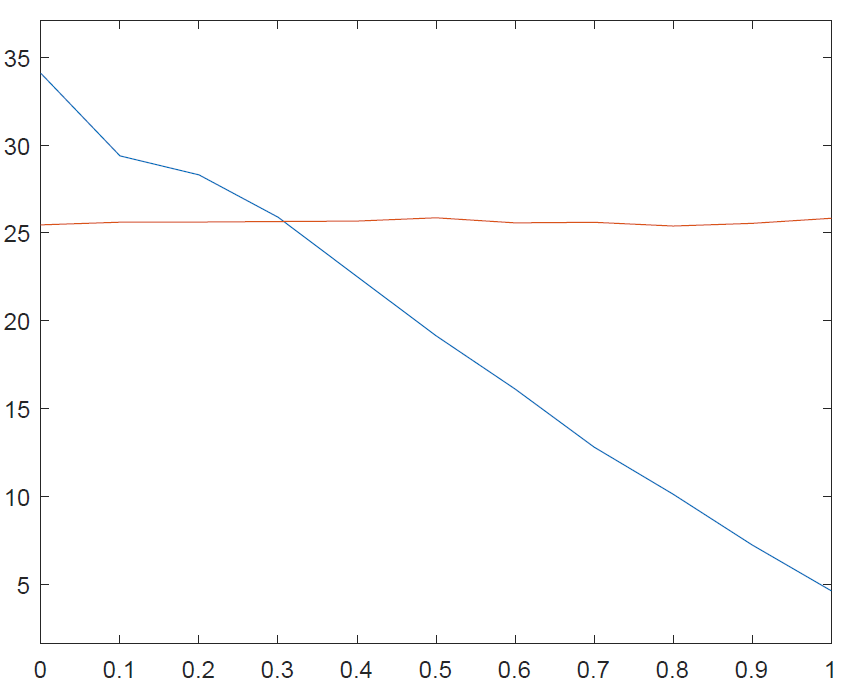}}
\caption{Зависимость меры различия \ref{def_difference_measure_of_degrees_distributions} от $\gamma$. Распределение степеней вершин. $\gamma_0=0.5$. Бразилия, 2012 год. 100 наблюдений. Красным - сеть вероятности совпадения знаков, синим - сеть Пирсона.}
\end{figure}
\begin{figure}[H]
\label{fig7}
\center{\includegraphics*[scale=0.7]{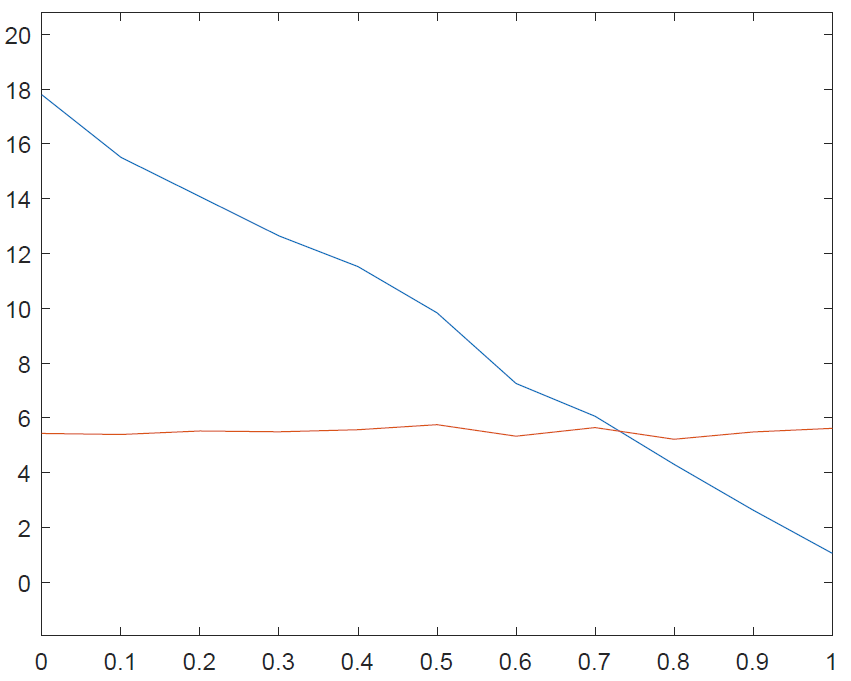}}
\caption{Зависимость меры различия \ref{def_difference_measure_of_degrees_distributions} от $\gamma$. Распределение степеней вершин. $\gamma_0=0.5$. Бразилия, 2012 год. 250 наблюдений. Красным - сеть вероятности совпадения знаков, синим - сеть Пирсона.}
\end{figure}
С другой стороны, есть примеры, когда сходимость обеих мер плохая и значение меры различия с увеличением наблюдений практически не изменяется.
\begin{figure}[H]
\label{fig8}
\center{\includegraphics*[scale=0.7]{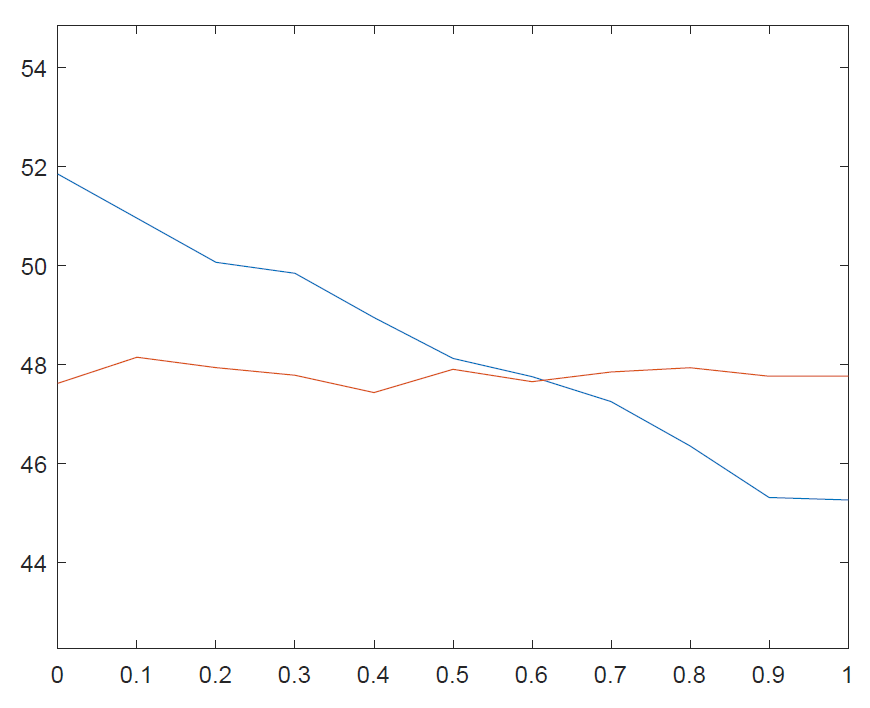}}
\caption{Зависимость меры различия \ref{def_difference_measure_of_degrees_distributions} от $\gamma$. Распределение степеней вершин. $\gamma_0=0.3$. Франция, 2011 год. 100 наблюдений. Красным - сеть вероятности совпадения знаков, синим - сеть Пирсона.}
\end{figure}
\begin{figure}[H]
\label{fig9}
\center{\includegraphics*[scale=0.7]{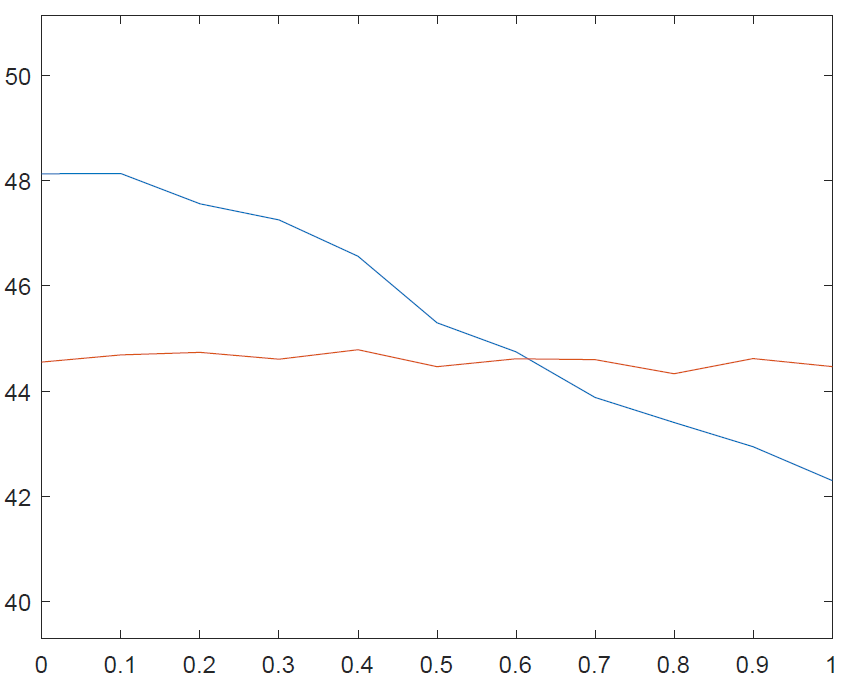}}
\caption{Зависимость меры различия \ref{def_difference_measure_of_degrees_distributions} от $\gamma$. Распределение степеней вершин. $\gamma_0=0.3$. Франция, 2011 год. 250 наблюдений. Красным - сеть вероятности совпадения знаков, синим - сеть Пирсона.}
\end{figure}
На рисунках ниже отмечены "экстремальные" случаи, когда мера различия при использовании для оценок вероятности совпадения знаков оказывается почти равной мере различия при использовании для оценок выборочного коэффициента корреляции Пирсона при $\gamma$ близких к $1$ и, наоборот, оказывается хуже при $\gamma$ близком к $0$.
\begin{figure}[H]
\label{fig10}
\center{\includegraphics*[scale=0.7]{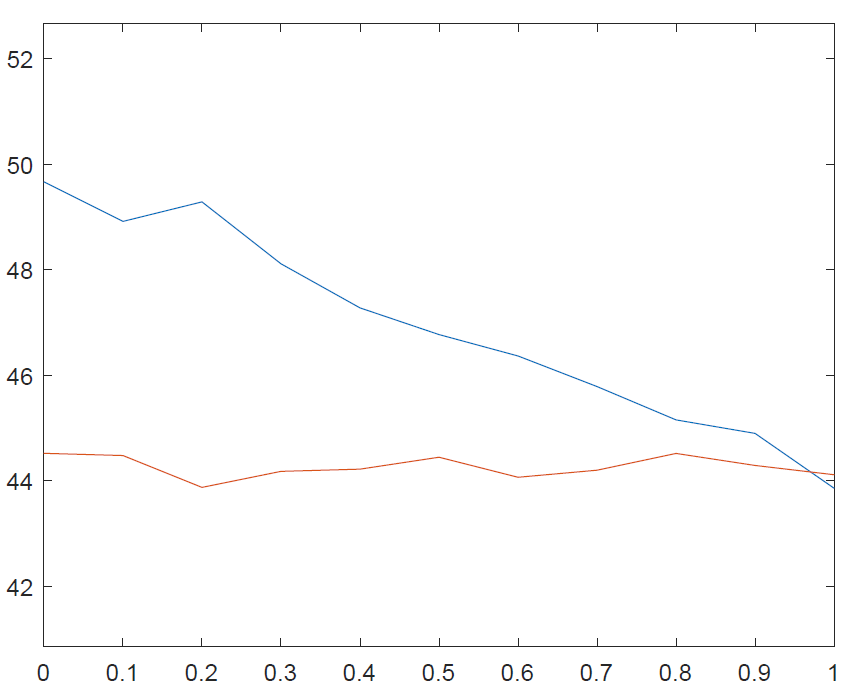}}
\caption{Зависимость меры различия \ref{def_difference_measure_of_degrees_distributions} от $\gamma$. Распределение степеней вершин. $\gamma_0=0.1$. Китай, 2011 год. 250 наблюдений. Красным - сеть вероятности совпадения знаков, синим - сеть Пирсона.}
\end{figure}
\begin{figure}[H]
\label{fig11}
\center{\includegraphics*[scale=0.7]{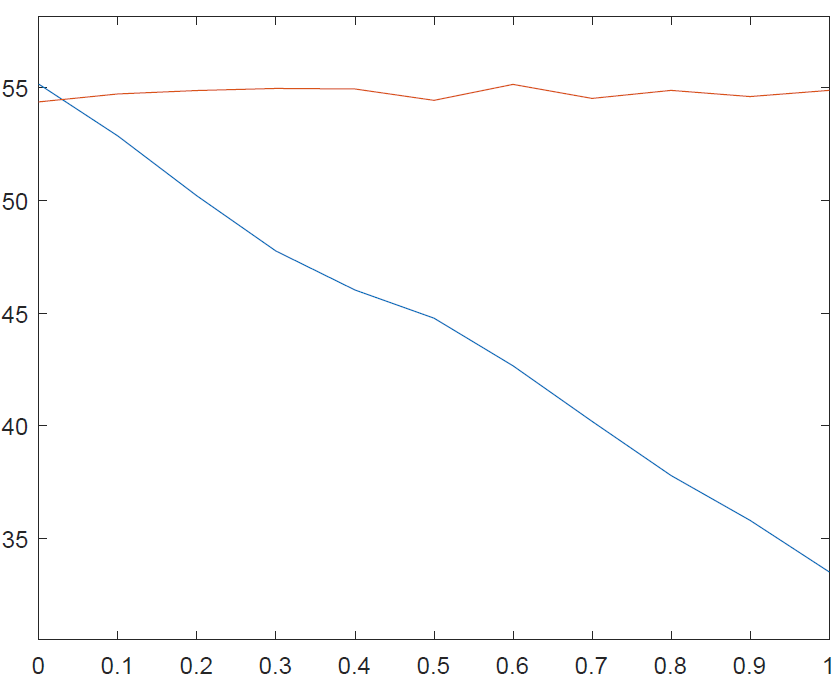}}
\caption{Зависимость меры различия \ref{def_difference_measure_of_degrees_distributions} от $\gamma$. Распределение степеней вершин. $\gamma_0=0.5$. Германия, 2010 год. 100 наблюдений. Красным - сеть вероятности совпадения знаков, синим - сеть Пирсона.}
\end{figure}
\subsection{Клики и независимые множества}
Эксперимент заключался в следующем: $n (n=100, 250)$  раз генерировался N-мерный $(N=50)$ случайный вектор с распределением смеси с заданным $\gamma$. По наблюдениям вычислялись выборочные корреляции Пирсона и частоты совпадения знаков, строился граф рынка \ref{def_threshold_graph} с порогом $\gamma_0 (\gamma_0=0,1; 0,3; 0,5)$, находилась максимальная клика (максимальное независимое множество). Эксперимент повторялся $10000$ раз и усреднением находилась оценка меры различия \ref{def_difference_measure_of_cliques}. Результаты проведения экспериментов представлены в виде кривых зависимости меры различия от $\gamma (\gamma=0, 0.1, 0.2, ... , 1)$  для каждого рынка, каждой сети и каждого значения числа наблюдений.
\\Для данных сетевых структур вероятность совпадения знаков оказывается устойчивой к изменению параметра смеси $\gamma$.
\begin{figure}[H]
\label{fig17}
\center{\includegraphics*[scale=0.7]{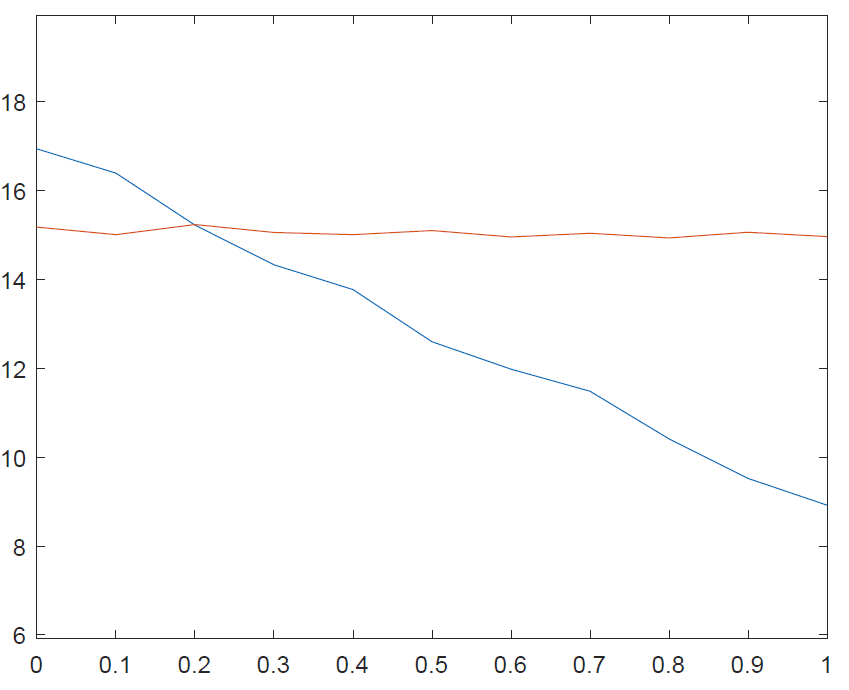}}
\caption{Зависимость меры различия \ref{def_difference_measure_of_cliques} от $\gamma$. Максимальное независимое множество. $\gamma_0=0.5$. Китай, 2011 год. 100 наблюдений. Красным - сеть вероятности совпадения знаков, синим - сеть Пирсона.}
\end{figure}
В данном случае характерно, что при увеличении числа наблюдений знаки показывают лучшую сходимость, нежели коэффициент корреляции Пирсона.
\begin{figure}[H]
\label{fig12}
\center{\includegraphics*[scale=0.7]{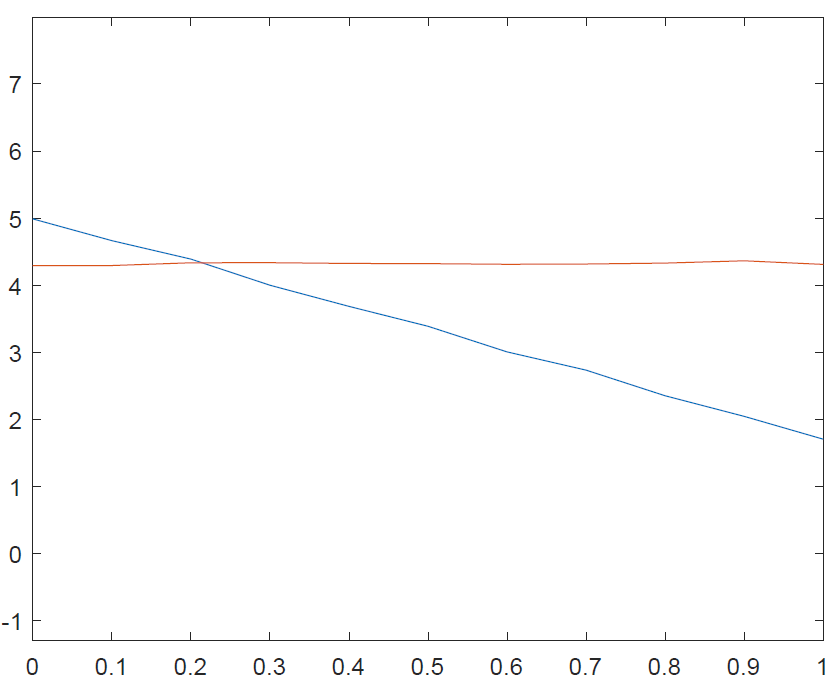}}
\caption{Зависимость меры различия \ref{def_difference_measure_of_cliques} от $\gamma$. Максимальная клика. $\gamma_0=0.1$. США, 2011 год. 100 наблюдений. Красным - сеть вероятности совпадения знаков, синим - сеть Пирсона.}
\end{figure}
\begin{figure}[H]
\label{fig13}
\center{\includegraphics*[scale=0.7]{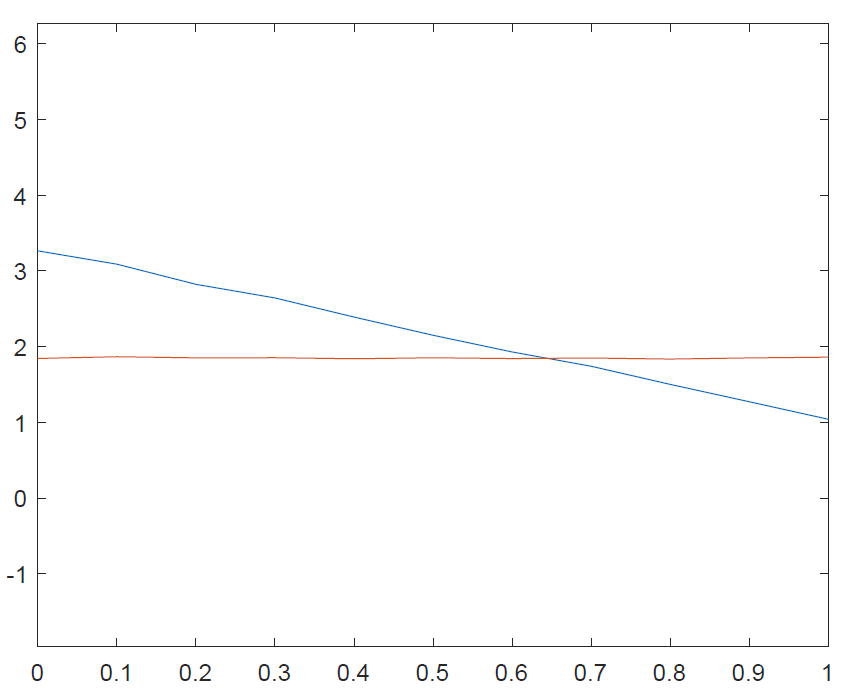}}
\caption{Зависимость меры различия \ref{def_difference_measure_of_cliques} от $\gamma$. Максимальная клика. $\gamma_0=0.1$. США, 2011 год. 250 наблюдений. Красным - сеть вероятности совпадения знаков, синим - сеть Пирсона.}
\end{figure}
И в этих экспериментах имеются примеры, когда мера различия при использовании для оценок вероятности совпадения знаков оказывается почти равной с мерой различия при использовании для оценок выборочного коэффициента корреляции Пирсона при $\gamma$ близких к $1$ и, наоборот, оказывается хуже при $\gamma$ близком к $0$.
\begin{figure}[H]
\label{fig14}
\center{\includegraphics*[scale=0.7]{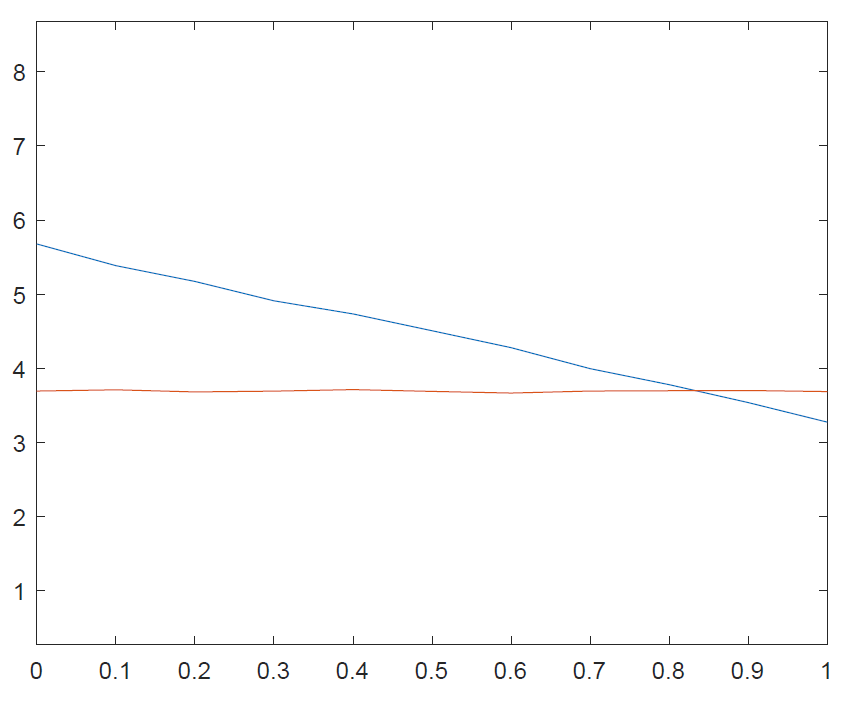}}
\caption{Зависимость меры различия \ref{def_difference_measure_of_cliques} от $\gamma$. Максимальная клика. $\gamma_0=0.5$. Индия, 2010 год. 250 наблюдений. Красным - сеть вероятности совпадения знаков, синим - сеть Пирсона.}
\end{figure}
\begin{figure}[H]
\label{fig15}
\center{\includegraphics*[scale=0.7]{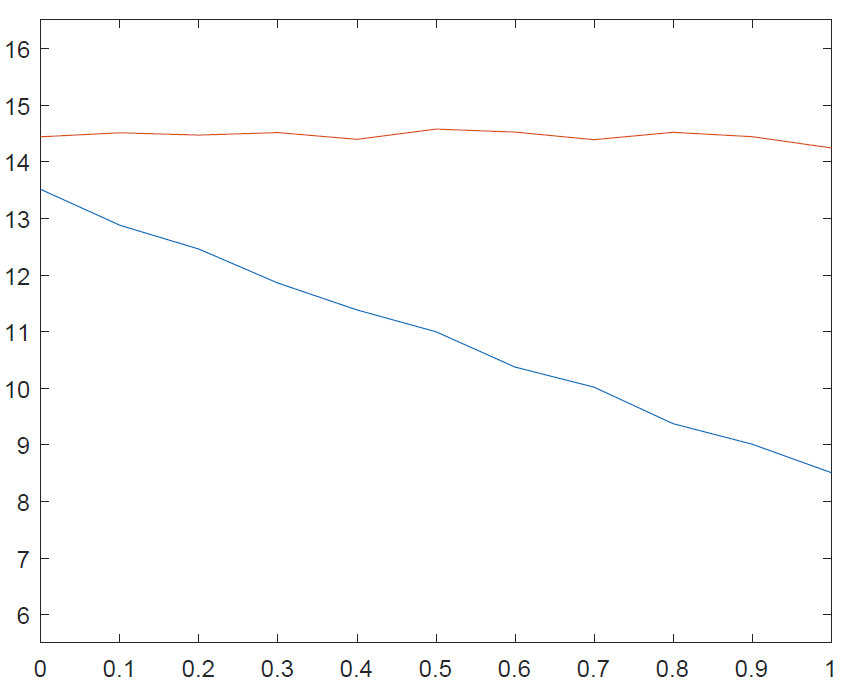}}
\caption{Зависимость меры различия \ref{def_difference_measure_of_cliques} от $\gamma$. Максимальная клика. $\gamma_0=0.1$. Россия, 2010 год. 100 наблюдений. Красным - сеть вероятности совпадения знаков, синим - сеть Пирсона.}
\end{figure}
\begin{figure}[H]
\label{fig16}
\center{\includegraphics*[scale=0.7]{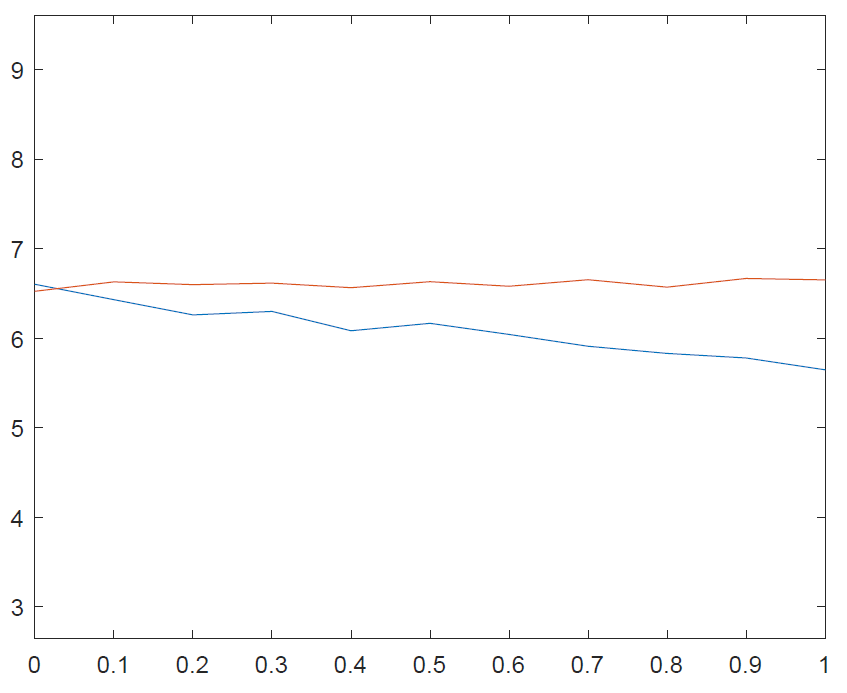}}
\caption{Зависимость меры различия \ref{def_difference_measure_of_cliques} от $\gamma$. Максимальное независимое множество. $\gamma_0=0.1$. Франция, 2011 год. 100 наблюдений. Красным - сеть вероятности совпадения знаков, синим - сеть Пирсона.}
\end{figure}
\subsection{Максимальное остовное дерево}
В данном случае рассматривалась смесь $\gamma*f_{gauss}(x_1,...,x_N) + (1-\gamma)*f_{student}(x_1,...,x_N)$ N-мерного нормального распределения и N-мерного распределения Стьюдента с 3 и 4 степенями свободы. Эксперимент заключался в следующем: $n (n=1000, 10000)$  раз генерировался N-мерный $(N=10)$ случайный вектор с распределением смеси с заданным $\gamma$. По наблюдениям вычислялись выборочные корреляции Пирсона и частоты совпадения знаков, находились MST \ref{def_MST}. Эксперимент повторялся $10000$ раз и усреднением находилась оценка меры различия \ref{def_difference_measure_of_MSTs}. Результаты проведения экспериментов представлены в виде кривых зависимости меры различия от $\gamma (\gamma=0, 0.1, 0.2, ... , 1)$  для каждого рынка, каждой сети, каждого значения количества наблюдений и каждой степени свободы распределения Стьюдента. Графики ниже показаны для смеси с распределением Стьюдента с 3 степенями свободы. Результаты для смеси с распределением Стьюдента с 4 степенями свободы аналогичны.
\\
При $n=1000$ значение меры различия (вероятности правильного определения топологии степеней вершин) очень близко к $0$.
\begin{figure}[H]
\label{fig21}
\center{\includegraphics*[scale=0.7]{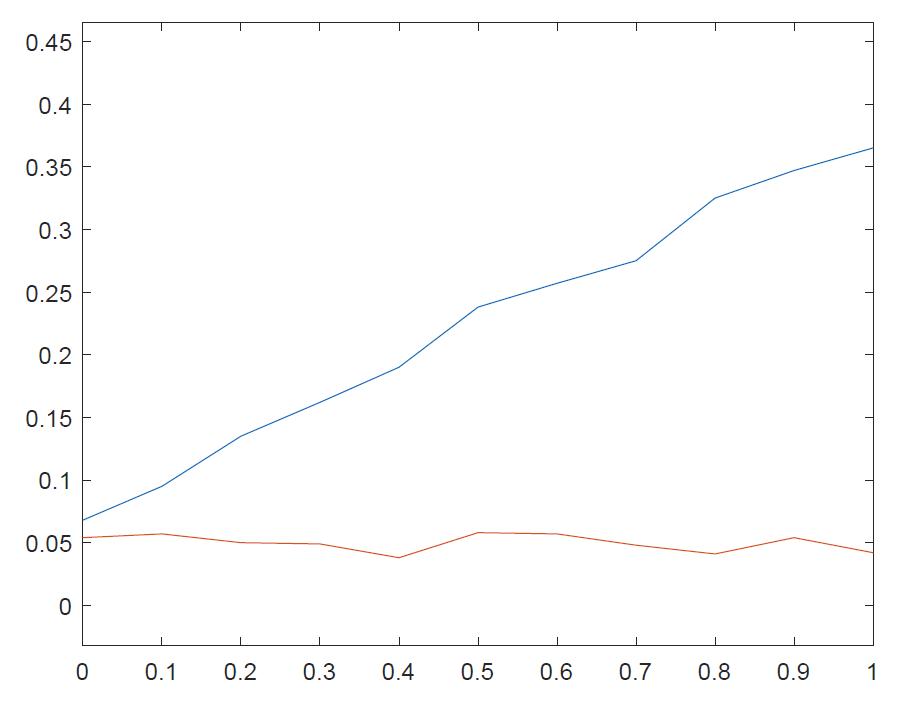}}
\caption{Зависимость меры различия \ref{def_difference_measure_of_MSTs} от $\gamma$ (смесь с распределением Стьюдента с 3 степенями свободы). Максимальное остовное дерево. Китай, 2012 год. 1000 наблюдений. Красным - сеть вероятности совпадения знаков, синим - сеть Пирсона.}
\end{figure}
При $n=10000$ графики для различных стран похожи друг на друга и основное различие заключается в значении меры различия. Наиболее характерные из графиков приведены ниже.
\begin{figure}[H]
\label{fig22}
\center{\includegraphics*[scale=0.7]{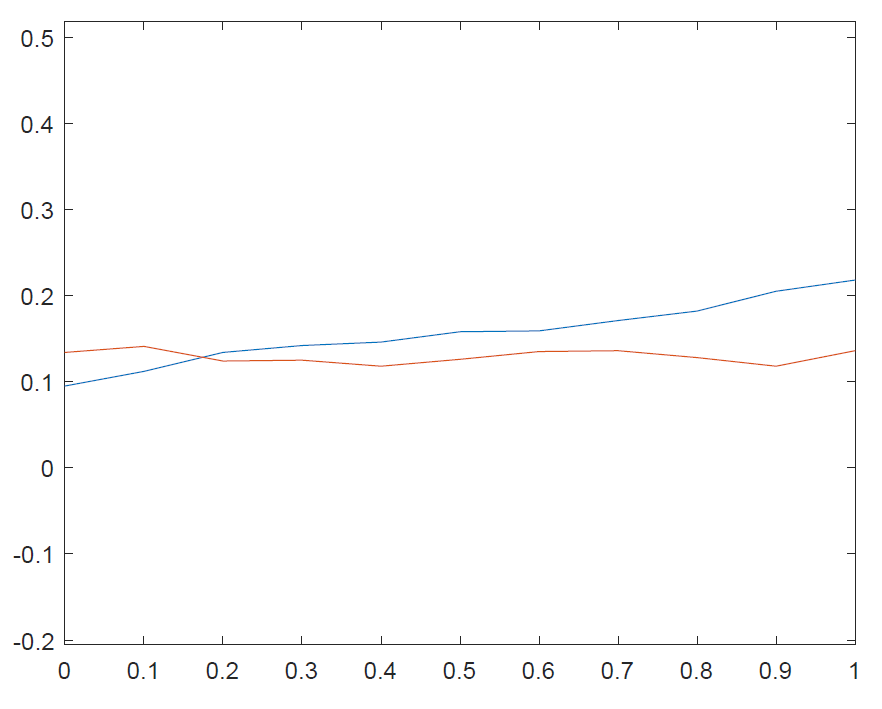}}
\caption{Зависимость меры различия \ref{def_difference_measure_of_MSTs} от $\gamma$ (смесь с распределением Стьюдента с 3 степенями свободы). Максимальное остовное дерево. Франция, 2010 год. 10000 наблюдений. Красным - сеть вероятности совпадения знаков, синим - сеть Пирсона.}
\end{figure}
\begin{figure}[H]
\label{fig23}
\center{\includegraphics*[scale=0.7]{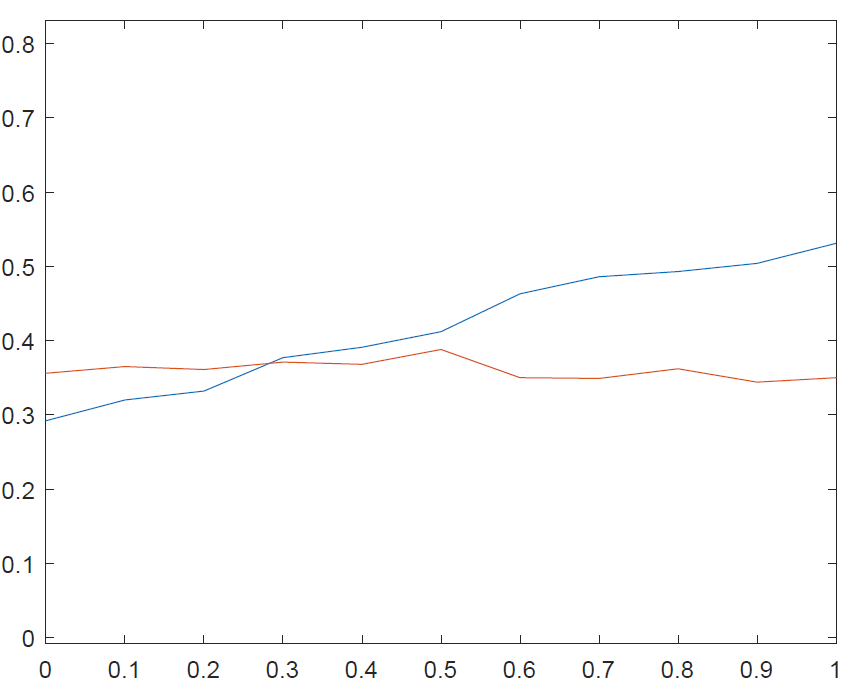}}
\caption{Зависимость меры различия \ref{def_difference_measure_of_MSTs} от $\gamma$ (смесь с распределением Стьюдента с 3 степенями свободы). Максимальное остовное дерево. Германия, 2012 год. 10000 наблюдений. Красным - сеть вероятности совпадения знаков, синим - сеть Пирсона.}
\end{figure}
\begin{figure}[H]
\label{fig24}
\center{\includegraphics*[scale=0.7]{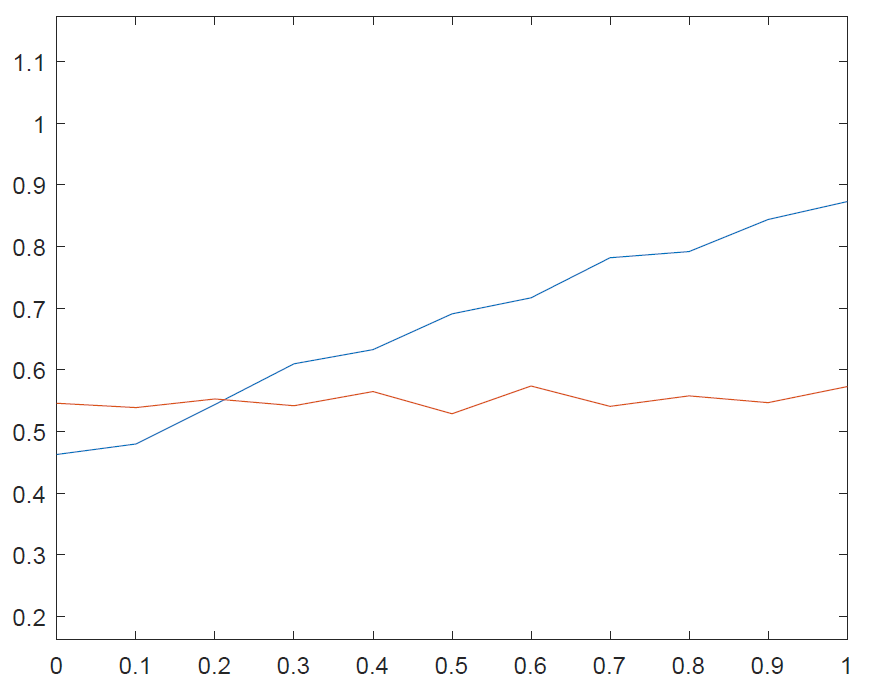}}
\caption{Зависимость меры различия \ref{def_difference_measure_of_MSTs} от $\gamma$ (смесь с распределением Стьюдента с 3 степенями свободы). Максимальное остовное дерево. Индия, 2011 год. 10000 наблюдений. Красным - сеть вероятности совпадения знаков, синим - сеть Пирсона.}
\end{figure}

Таким образом, приведенные результаты показывают устойчивость знаковых процедур оценивания характеристик сетевых моделей и неустойчивость процедур, основанных на выборочных корреляциях Пирсона, к изменению вероятностной модели совместного распределения доходностей акций в классе эллиптических распределений.
Работа выполнена при финансовой поддержке Российского гуманитарного научного фонда (проект 15-32-01052).

\end{document}